# A homotopy theory approach to set theory

notes by misha gavrilovich

> *but the mischief of it is, nature will have to take its course*
>
> *Miguel de Cervantes Saavedra, Don Quixote*

**Abstract**:

We observe that the notion of *two sets being equal up to finitely many elements* is a *homotopy equivalence* relation in a model category, a common axiomatic formalism for homotopy theory introduced by Quillen "to cover in a uniform way a large number of arguments in homotopy theories that were formally similar to well-known ones in algebraic topology. We show the same formalism covers some arguments in (naive) set theory, and a well-known set-theoretic invariant, the covering number $\mathrm{cf}([\aleph_\omega]^{\aleph_0})$, of PCF theory. Further we observe a similarity between homotopy theory ideology/yoga and that of PCF theory, and briefly discuss conjectural connections with model theory and arithmetics and geometry.

We argue that the formalism is curious as it suggests to look at a *homotopy-invariant variant of Generalised Continuum Hypothesis* about which more can be proven within ZFC and first appeared in PCF theory independently but with a similar motivation.

1. This is a dense announcement of results partly reporting on joint work with Assaf Hasson, and shall eventually appear in the form of a joint paper. We give no proofs but the definitions are given in full detail. Proofs, speculations and motivations may be found in a more verbose report [Gavrilovich] which also contains some questions.

1.1. **Introduction.** *The structure of the paper.* In §1 and §2 we develop a notion of homotopy between sets such that two sets are *homotopic* iff they are *almost* equal, i.e. if they differ by finitely many elements. We do so in the formalism of a Quillen model category. In §1.3.2 we define a model category $QtNaamen$ of (some) families of sets, such that two objects $\{A\}$ and $\{B\}$ are *weakly (homotopy) equivalent* iff $A$ and $B$ are almost equal. In §1.2 we informally explain how to obtain the definition of $QtNaamen$ by starting from the poset(class) of all sets ordered by inclusion (considered as a category): declare some inclusions $A \subseteq B$ as weak equivalences and/or cofibrations so as to capture notions of finiteness, countability and equicardinality, and then formally add morphisms that have necessary lifting properties. In §3 we define and calculate the derived functor of the function $\operatorname{card} : QtNaamen \longrightarrow On$ sending a family of sets into its cardinality, and identify its values with a set-theoretic invariant in PCF theory, the covering number $\operatorname{cov}(\kappa, \aleph_1, \aleph_1, 2)$. In §1.3.3 and §2.2-2.7 we conjecture a similar construction gives a model category associated to a (not necessarily first order) uncountably categorical theory, and discuss algebro-geometric applications. In §3.8 we make some remarks and indicate directions for further research. All necessary definitions, both of homotopy and set theory, are given; appendices provide some examples as well as a definition of a model category. We work in naive set theory and expressedly ignore the hindrance of a particular axiomatisation of set theory; a reader may assume that we work either in a Grothendieck universe, a model of ZFC or a set theory with a universal set, cf. §3.8.3.

1.2. *A sketch of the construction of the model category.* A reader not familiar with model categories may skim through the explanation at the first reading.

1.2.1. We start by considering, as a category, the partially ordered set(class) of all sets ordered by inclusion. We want an inclusion $A \subseteq B$ to be a weak (homotopy) equivalence iff $B \setminus A$ is finite. A model category is a category equipped with the extra structure consisting of the three classes of morphisms called(labelled as) (*w*) weak (homotopy) equivalences, (*c*) cofibrations and (*f*) fibrations. Further we define an inclusion $A \subseteq B$ to be both a weak (homotopy) equivalence and a cofibration iff $B \setminus A$ is finite; this is motivated by symmetry considerations as well as as a desire to follow the standard construction of a cofibrantly generated model structure.

1.2.2. However, we lack fibrations, although e.g. Axiom M2 implies that every morphism is a fibration up to a weak equivalence but no morphism-inclusion $A \subseteq B$ right-lifts wrt to $(wc)$-inclusion $A \subseteq A \cup \{b\}$ for $b \in B \setminus A$ as required by Axiom M1 of a fibration. Therefore we introduce the fibrations formally, as formal morphisms having the required lifting properties, in the following way.

1.2.3. By Axiom M2, every morphism-inclusion $A \subseteq B$, write $A \longrightarrow B$, decomposes as a composition $A \xrightarrow{(wc)} A' \xrightarrow{(f)} B$ where the unique morphism $A \longrightarrow A'$ is both a (w) weak equivalence as well as a (c) cofibration, and the unique morphism $A' \longrightarrow B$ is a (f) fibration. We add $A' \longrightarrow B$ as a new (for $A \neq B$) morphism that by definition satisfies the necessary lifting properties, particularly it left-lifts wrt to $A \longrightarrow A \cup \bar{b}$ for all finite $\bar{b} \subseteq B$. We also add $A \longrightarrow A'$ as a new (for $A \neq B$) morphism.

1.2.4. As (we insist that) there is at most one morphism between any pair of objects, adding formally a morphism $A' \longrightarrow B$ is the same as adding a formal object $A'$, and the latter can be identified with a *family* of subsets of $B$ (of the form $A \cup \bar{b}, \bar{b} \subseteq B$ finite).

1.2.5. To enable our model category to express *(equi)cardinality* and *countability*, we declare an inclusion $A \subseteq B$ a cofibration iff either $\operatorname{card} A = \operatorname{card} B$ or both $A$ and $B$ are at most countable. Similarly to the above, we then add formally $\xrightarrow{(wf)}$-morphisms (i.e. those that are both a fibration and a weak equivalence).

1.2.6. Thus, to get a model category we need to add (at least some) families of sets as formal objects. In the construction of *StNaamen* below, we add *all* families of sets as formal objects; in the construction of the model category *QtNaamen* which is a full subcategory of *StNaamen*, we are more careful and add formally only some families of sets. Sets $A$ and $B$ therefore become singleton families $\{A\}$ and $\{B\}$, and we write $\{A\} \longrightarrow \{B\}$ for the inclusion $A \subseteq B$.

1.3. **Definition** *(c-w-f arrow notation and the construction of the model category).* We write $A \xrightarrow{g} B \measuredangle X \xrightarrow{h} Y$ iff for any $A \longrightarrow X$ and $B \longrightarrow Y$ there exists $B \longrightarrow X$ and $A \longrightarrow Y$ such that

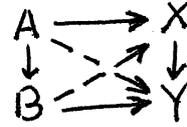

(⊠) if the square of solid arrows commutes, then the whole diagram commutes.

We read $A \longrightarrow B \measuredangle X \longrightarrow Y$ as: the morphism (or arrow) $A \longrightarrow B$ lifts wrt $X \longrightarrow Y$. We use dotted arrows to indicate an existential quantifier. In the categories we shall consider there is at most one morphism between any two objects, and therefore any diagram that can drawn is necessarily commutative, e.g. we would not need to check the condition (⊠) above.

1.3.1. Let *StNaamen* be the following labelled category. Its objects are arbitrary sets, and there is at most one arrow between any two objects. An arrow carries none or some of the three labels c, w, f, and we write e.g. $A \xrightarrow{(wc)} B$ to indicate that the arrow $A \longrightarrow B$ carries labels w, c and possibly f. Auxiliary items $(\rightarrow)_0$, $(wc)_0$, $(c)_0$ define notation used to bootstrap the definition, and labels $(wc)_0$, $(c)_0$ are not part of structure of *StNaamen*. We define:

$(\rightarrow)_0$ $\{A\} \longrightarrow \{B\}$ iff $A \subseteq B$

$(wc)_0$ $\{A\} \xrightarrow{(wc)_0} \{B\}$ iff the difference $B \setminus A$ is finite (and $A \subseteq B$)

$(c)_0$ $\{A\} \xrightarrow{(c)_0} \{B\}$ iff card $A$ = card $B$ or card $B \leqslant \aleph_0$ (and $A \subseteq B$)

$(\rightarrow)$ $X \longrightarrow Y$ iff $\forall x \in X \ \exists y \in Y \ (x \subseteq y)$

$(f)$ $X \xrightarrow{(f)} Y$ iff $\{a\} \longrightarrow \{b\} \measuredangle X \longrightarrow Y$ whenever $\{a\} \xrightarrow{(wc)_0} \{b\}$

$(wf)$ $X \xrightarrow{(wf)} Y$ iff $\{a\} \longrightarrow \{b\} \measuredangle X \longrightarrow Y$ whenever $\{a\} \xrightarrow{(c)_0} \{b\}$

$(c)$ $A \xrightarrow{(c)} B$ iff $A \longrightarrow B \measuredangle X \longrightarrow Y$ whenever $X \xrightarrow{(wf)} Y$

$(wc)$ $A \xrightarrow{(wc)} B$ iff $A \longrightarrow B \measuredangle X \longrightarrow Y$ whenever $X \xrightarrow{(f)} Y$

$(w)$ $A \xrightarrow{(w)} Y$ iff it decomposes as $A \xdashrightarrow{(wc)} \cdot \xdashrightarrow{(wf)} Y$

©2010

1.3.2. Let $QtNaamen \subseteq StNaamen$ be the full subcategory of $StNaamen$ consisting of sets $X$ such that any of the following equivalent conditions hold for any $M, a, b$ and arrows $A \xrightarrow{(c)} B$, $B' \xrightarrow{(wf)} B$, $X \longrightarrow Y$:

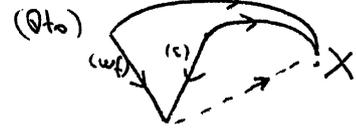

$(Qt_0)$

$(Qt_1)$ $A \xrightarrow{(c)} B \measuredangle X \longrightarrow Y$ or $B' \xrightarrow{(wf)} B \measuredangle X \longrightarrow Y$

$(Qt'_1)$ $\cup \{X'' : X \dashleftarrow X_0 \xleftarrow{(c)} X'' \xrightarrow{(wf)} X' \dashrightarrow X\} \dashrightarrow X$

$(Qt_2)$ if $\{a\} \to X$, $B' \longrightarrow X$ and $\{a\} \xrightarrow{(c)} \{b\}$ and $B' \xrightarrow{(wf)} \{b\}$, then $\{b\} \longrightarrow X$

$(Qt_3)$ $M^{\leq \aleph_0} \longrightarrow X$ implies $M^{\leq \mathrm{card}(x \cap M)} \longrightarrow X$ for any $x \in X$ (where $M^{\leq \lambda} := \{L \subseteq M : \mathrm{card}\, L \leq \lambda\}$)

1.3.3. Let $\mathfrak{R}_\preceq = (\mathfrak{K}, \preceq_\mathfrak{R})$ be a class of models equipped with a partial order relation $\preceq_\mathfrak{R}$, e.g. $\mathfrak{R} = \mathfrak{R}(T)$ be the class of models of a theory $T$ (in a fixed language $L = L(T)$), and for models $M, N \in \mathfrak{R}(T)$ $M \preceq_\mathfrak{R} N$ iff $M \subseteq N$ is an elementary submodel of $N$. Let $QtNaamen(\mathfrak{R}_\preceq)$ be the subcategory of $QtNaamen$ of families of models in $\mathfrak{R}$ with a fibrantly generated model structure. That is, it is defined as follows:

$\mathcal{O}b(\mathfrak{R}_\preceq)$ $X \in \mathcal{O}bQtNaamen(\mathfrak{R}_\preceq)$ iff $\forall M \in X \ (M \in \mathfrak{R})$

$(\to_{\preceq_\mathfrak{R}})$ $X \longrightarrow Y$ iff $\forall M \in X \ \exists N \in Y \ (M \preceq N)$

$(f)$ $X \xrightarrow{(f)} Y$ iff $X \xrightarrow{(f)Qt} Y$ is an (f)-arrow in QtNaamen

$(wf)$ $X \xrightarrow{(wf)} Y$ iff $X \xrightarrow{(wf)Qt} Y$ is an (wf)-arrow in QtNaamen

$(c)$ $A \xrightarrow{(c)} B$ iff $A \longrightarrow B \measuredangle X \longrightarrow Y$ whenever $X \xrightarrow{(wf)} Y$ (in $QtNaamen(\mathfrak{R}_\preceq)$)

$(wc)$ $A \xrightarrow{(wc)} B$ iff $A \longrightarrow B \measuredangle X \longrightarrow Y$ whenever $X \xrightarrow{(f)} Y$ (in $QtNaamen(\mathfrak{R}_\preceq)$)

$(w)$ $A \xrightarrow{(w)} Y$ iff it decomposes as $A \xdashrightarrow{(wc)} \cdot \xdashrightarrow{(wf)} Y$ (in $QtNaamen(\mathfrak{R}_\preceq)$)

2. **Model categories. A theorem. A conjecture.** The above c-w-f-arrow notation allows one to use in set theory, and possibly model theory, the language of commutative diagrams of model categories, e.g. to draw an analogy between a fibre bundle $V \longrightarrow B$ and an inductive construction $\{M_i\}_i \xrightarrow{(f)} \{\cup_i M_i\}$, or view Lowenheim-Skolem theorem as an instance of an $\cdot \xrightarrow{(c)} \cdot \xrightarrow{(wf)} \cdot$-decomposition required by Axiom M2. Here $f$ stands for fibration, $c$ stands for cofibration, and $w$ stands for weak (homotopy) equivalence.

Below we develop an example to show that the language of model categories retains some of its powers. Diagram chasing along with basic set theory arguments gives:

2.1. **Theorem.** *The category QtNaamen with these labels is a model category* (see the appendix for the definition).

2.2. Further we conjecture that often a (not necessarily first-order) *categorical* theory, via the class of (families of) its models, defines a model subcategory of $QtNaamen$

2.3. **The Eyjafjallajoekull conjecture** (Bays and Gavrilovich). *Let $\mathfrak{R}_{\preceq} = (\mathfrak{K}, \preceq_{\mathfrak{R}})$ be an excellent abstract elementary class. Then there is a subcategory $Qt(\mathfrak{R}_{\preceq})$ of $QtNaamen(\mathfrak{R}_{\preceq})$ which is a model category, and is not degenerate, e.g. for every cardinality $\lambda$ there is a model $M_\lambda \in \mathfrak{R}$ of cardinality $\lambda$ such that $\{M_\lambda\} \in \mathcal{O}bQt(\mathfrak{R}_{\preceq})$, and for any $M, N \in \mathfrak{R}$ it holds that $\{M\} \xrightarrow{(wc)} \{N\}$ iff $M \preceq N$ and $N$ is primary over $M \cup \bar{b}$ for a finite set $\bar{b} \subset N$.*

2.4. Axiom M5(2-out-of-3) is delicate to check and seem to impose structural constrains; mere diagram chasing suffices to prove that *any complete quasi-poset with fibrantly or cofibrantly generated model category labelling, and particularly the categories $StNaamen(\mathfrak{R}_{\preceq})$ and $QtNaamen(\mathfrak{R}_{\preceq})$ with labels as defined above, satisfy Axioms M1-M4 and M6 of model categories* but not necessarily M5(2-out-of-3).

2.5. Our intention is that the definition of $Qt(\mathfrak{R}_{\preceq})$ is to be such that model category diagram chasing corresponds to arguments in excellent classes, e.g. drawing the pushout square on the left is to correspond to a categoricity transfer argument in excellent classes (§3,p.17 of [Lessmann]). Take a model $M$ of cardinality $\aleph_1$, take its cofibrant replacement $\emptyset \xrightarrow{(c)} \{M_i\}_{i\in\omega_1} \xrightarrow{(wf)} \{M\}$ by splitting $M = \cup_{i\in\omega_1} M_i$ into a continuous increasing chain of countable models $M_i$. Pick an element $a$ and construct an acyclic cofibration $\{M_i\}_i \xrightarrow{(wc)} \{\overline{M_i a}\}_i$ of models $\overline{M_i a}$ countable primary over $M_i \cup \{a\}$. Finally, take the pushout of $\{M\} \xleftarrow{(wf)} \{M_i\}_{i\in\omega_1} \xrightarrow{(wc)} \{\overline{M_i a}\}_{i\in\omega_1}$. In $StNaamen(\mathfrak{R}_{\preceq})$ the pushout is simply $\{M\} \cup \{\overline{M_i a}\}_{i\in\omega_1}$, whereas in $QtNaamen(\mathfrak{R}_{\preceq})$ the pushout is a single model $\{N\}$, $N \supseteq M \cup \{a\}$. By Axioms M3 and M4 or M5 we get $\{M\} \xrightarrow{(wc)} \{N\}$, and thus $N$ is primary over $M \cup \{a\}$.

2.6. The conjecture appears to relate model categories and such questions as *Mumford-Tate, Kummer theory, Mordell-Weil, Schanuel conjecture and dual thereto*, as these algebro-geometric properties are essentially used to prove excellency and uncountable categoricity (rather, finitely many models) of the following explicitly given abstract elementary classes of models, cf. [Bays, DPhil] and references therein.

2.7. Fix an elliptic curve $E/\bar{\mathbb{Q}}$ without complex multiplication defined over a number field. Consider the classes of all short exact sequences of the form

$$(\overline{K}^*) \quad 0 \longrightarrow \mathbb{Z} \longrightarrow V \xrightarrow{\varphi} \overline{K}^* \longrightarrow 1$$
$$(\overline{K}_p^*) \quad 0 \longrightarrow \mathbb{Z}[\tfrac{1}{p}] \longrightarrow V \xrightarrow{\psi} \overline{K}_p^* \longrightarrow 1$$
$$(E) \quad 0 \longrightarrow \mathbb{Z}^2 \longrightarrow V \xrightarrow{\varphi} E(\overline{K}) \longrightarrow 0$$

©2010

($\mathbb{C}_{\exp}$) a field $\overline{K}$ equipped with a homomorphism $\exp : \overline{K} \longrightarrow \overline{K}^*$

where, as notation suggests, $V$ varies among $\mathbb{Q}$-vector spaces, $\overline{K}$ varies among algebraically closed fields of zero characteristic, and $\varphi$ varies among group homomorphisms; and $\overline{K}_p$ varies among algebraically closed fields of prime characteristic $p$, and $\psi$ varies among $\mathbb{Z}[\frac{1}{p}]$-module homomorphisms such that the isomorphism type of the restriction $\psi_{|\mathbb{Q}\mathbb{Z}[\frac{1}{p}]} : \mathbb{Q}\mathbb{Z}[\frac{1}{p}] \longrightarrow \bar{\mathbb{F}}_p$ is fixed. For the first three classes, let $\preceq_\Re$ be all inclusions of submodels respecting field and vector space structure; these are necessarily elementary. The definition of class ($\mathbb{C}_{\exp}$) is more complicated and may be found elsewhere [Zilber, Pseudoexponentiation]. All four are excellent abstract elementary classes [ibid.,Bays] that have finitely many, up to isomorphism, models of each uncountable cardinality, i.e. for every uncountable $\bar{K}$ there are finitely many short exact sequences up to a linear isomorphism of $V$ inducing a field automorphism on $\bar{K}$.

3. **An application: PCF theory as a homotopy-invariant theory.** We derive the function $\mathrm{card} : QtNaamen \longrightarrow On$, observe that its homotopy-invariant left derived functor $\mathbb{L}_c \mathrm{card} : QtNaamen \longrightarrow On$ is studied in PCF theory under the name of *the covering number*, and make some remarks on similarity between ideology/yoga of homotopy theory and PCF theory.

3.1. Let $\mathcal{A}, \mathcal{B}$ be quasi-partially ordered sets considered as categories where $x \longrightarrow y$ iff $x \leqslant y$. Then a (covariant) *functor* $F : \mathcal{A} \longrightarrow \mathcal{B}$ is a non-decreasing function $F : \mathcal{A} \longrightarrow \mathcal{B}$. If both $\mathcal{A}$ and $\mathcal{B}$ are also equipped with a c-w-f labelling, we say that a functor $F : \mathcal{A} \longrightarrow \mathcal{B}$ is *homotopy-invariant* iff for any arrow $X \xrightarrow{(w)} Y$ (weak homotopy equivalence), it holds $F(X) \xrightarrow{(w)} F(Y)$. An initial object $\perp$ of $\mathcal{A}$ is a minimal element of $\mathcal{A}$ (whenever such exists). (As any diagram is commutative in these categories, we need not state the conditions that the functors have to respect commutative diagrams.)

3.2. Let $On$ be the category of ordinals where each arrow is labelled $(cf)$ and each isomorphism is labelled $(cwf)$. For a *function* $F : \mathcal{A} \longrightarrow On$, define (minimum is taken over all finite sequences labelled as shown)

$$\mathbb{L}_c F(X) = \min \left\{ F(Y) : \begin{array}{c} X_1 \\ \nearrow \ \nwarrow_{(w)} \end{array} \begin{array}{c} X_3 \\ \nearrow \ \nwarrow_{(w)} \end{array} \begin{array}{c} X_n \dashrightarrow Y \\ \nearrow \quad \uparrow \\ \quad | (c) \\ \quad | \\ \perp \end{array} \right\}$$

3.3. $\mathbb{L}_c F(X)$ is a homotopy invariant functor "closest from the left"(Quillen, I:4.1) to the function $F : StNaamen \longrightarrow On$, by which is meant: for any homotopy-invariant functor $G : StNaamen \longrightarrow On$ such that $G(X) \longrightarrow F(X)$ for any object $X$ such that $\perp \xrightarrow{(c)} X$, it holds that $G(Y) \longrightarrow \mathbb{L}_c F(Y)$ for any $\perp \xrightarrow{(c)} Y$ (note then there is a natural transformation from functor $G$ to functor $\mathbb{L}_c F$).

In particular, the function $\mathbb{L}_c F : StNaamen \longrightarrow On$ is the left derived functor of $F : StNaamen \longrightarrow On$ provided that $F$ is a functor.

© 2010

3.4. Take $F = \text{card}$ to be the cardinality function. Arguably, the model category formalism suggests we view $\mathbb{L}_c\text{card} : StNaamen \longrightarrow On$ as an analogue of a cofibrantly replaced left derived functor of the "forgetful functor" $\text{card} : StNaamen \longrightarrow On$. Then homotopy yoga suggests we view values of $\mathbb{L}_c\text{card}$, e.g. $\mathbb{L}_c\text{card}(\{\aleph_\alpha\}) = \mathbb{L}_c\text{card}(\{X : X \subseteq \aleph_\alpha\})$, as (homotopy-invariant and therefore) more robust and interesting invariants, as compared with the non-homotopy-invariant values $\text{card}(\{X : X \subseteq \aleph_\alpha\})$.

3.5. And indeed, it is for the reasons of being more robust and less prone to change by forcing that the values of $\mathbb{L}_c\text{card}(\{\aleph_\alpha\})$ (for limit $\aleph_\alpha$) have been introduced in set theory (Shelah, Cardinal Arithmetic). Set-theoretically, $\mathbb{L}_c\text{card}(\{\aleph_\alpha\}) = \text{cov}(\aleph_\alpha, \aleph_1, \aleph_1, 2)$ is the least size of a family $X$ of countable subsets of $\aleph_\alpha$, such that every countable subset of $\aleph_\alpha$ is a subset of a set in the family $X$. This may used, for example, to study the cardinality $(\aleph_\alpha)^{\aleph_0}$ of the set of countable subsets of $\aleph_\alpha$, via the bound $(\aleph_\alpha)^{\aleph_0} \leqslant \text{cov}(\aleph_\alpha, \aleph_1, \aleph_1, 2) + 2^{\aleph_0}$, by decomposing it into a "noise" "non-homotopy-invariant" part $2^{\aleph_0}$ whose value is known to be highly independent of ZFC (and easy to force to change), and a homotopy-invariant part $\text{cov}(\aleph_\alpha, \aleph_1, \aleph_1, 2)$ which admit bounds in ZFC (and is harder to force to change).

3.6. A short calculation gives $\mathbb{L}_c\text{card}(\{X : X \subseteq \aleph_0\}) = 1$ (in ZFC) whereas it is known that there are models of ZFC where e.g. $\text{card}(\{X : X \subseteq \aleph_0\}) = 2^{\aleph_0} > \aleph_{\omega_\omega}$. Meanwhile, non-trivially, Shelah (Cardinal Arithmetic, IX:4) proves $\mathbb{L}_c\text{card}(\{\aleph_\omega\}) < \aleph_{\omega_4}$. Similar upper bounds exist on $\mathbb{L}_c\text{card}(\{\aleph_\alpha\})$ for (most) $\aleph_\alpha$ limit (excepting $\aleph_\alpha = \alpha$), and are provided by PCF theory.

3.7. Arguably, the above justifies saying that the homotopy-invariant version of Generalised Continuum Hypothesis has less independence of ZFC, as suggested by homotopy theory.

3.8. **Remarks.** These remarks are explained in more details in [Gavrilovich].

3.8.1. Gromov [Ergosystems] writes that "The category/functor modulated structures can not be directly used by ergosystems, e.g. because the morphisms sets between even moderate objects are usually unlistable. But the ideas of the category theory show that there are certain (often non-obviuos) rules for generating proper concepts." Curiously, in our categories where this obstruction does not arise, all definitions we make seem to be a result of a rather direct and automatic, straightforward repeated application of the lifting property to basic concepts of naive set theory, and the axioms of a model category admit a functional semantics whereby they are interpreted as rules to draw arrows and add labels on labelled graphs. We say more on this in [Gavrilovich], particularly §1.0.4,p.5 and §1.3,pp.12-14.

3.8.2. Shelah explicitly states his ideology of PCF theory in Shelah (Logical Dreams), e.g. Thesis 5.10, and we find it remarkably similar to the model category ideology as applied to StNaamen. It is unclear whether a deeper connection with PCF theory exists, e.g. whether the sequence of PCF generators is a (non-pointed, non-functorial) analogue of a (co)fibration sequence, or whether $X \longmapsto \{X\}$ and $X \longmapsto \cup_{x \in X} x$ can be usefully viewed as analogues of suspension $X \longmapsto \Sigma X$ and loop $X \longmapsto \Omega X$ spaces, cf. Kojman (A short proof of PCF theorem).

3.8.3. Manin (A course in logic, 2010, p.174) discusses the Continuum Hypothesis and the possibility for a need to "try to find alternative languages and semantics" for set theory. It would seem that the connection between homotopy theory (in the model category formalism) and set theory (in ZFC or NF, or similar formalisms) we suggest, may provide for such an alternative language and semantics.

3.8.4. Note that a topological space $T$ determines a homotopy-invariant functor $\mathsf{acc}_T : QtNaamen \longrightarrow Naamen$, $\mathcal{X} \longmapsto \cup_{X \in \mathcal{X}} \mathsf{acc}_T(X \cap T)$ sending a family $\mathcal{X}$ into the set of accumulation points $\cup_{X \in \mathcal{X}} \mathsf{acc}_T(X \cap T)$ of a member of the family; here $Naamen$ is the poset of all sets under inclusion. It appears that the definition of a topological space may be stated purely category-theoretically in terms of this functor and the functor $\{\cdot\} : Naamen \longrightarrow QtNaamen, X \longmapsto \{X\}$.

3.8.5. Our original motivation was to associate a model category (via the class of families of models) to an uncountably categorical theory and, more generally, to an excellent abstract elementary class (Shelah, Classification theory of non-elementary classes). In particular, we wanted to use the language of homotopy theory to perform the model-theoretic analysis of complex exponentiation $(\mathbb{C}, +, *, \exp)$ (Zilber, Pseudo-exponentiation on algebraically closed fields of characteristic zero) and covers of semi-Abelian varieties ([Bays] and references therein). These results claim there exist a unique, up to an appropriate notion of isomorphism (*not* respecting topology), function $\mathrm{ex} : \mathbb{C} \longrightarrow \mathbb{C}$ satisfying $\mathrm{ex}(x+y) = \mathrm{ex}(x)\mathrm{ex}(y)$, the Schanuel conjecture and a dual thereto; Bays replaces $\mathbb{C}$ and $\mathrm{ex}$ by an elliptic curve and its cover $\mathrm{ex}_E : \mathbb{C} \longrightarrow \mathbb{C}/\Lambda$. Their analysis leads to a number- and geometric-theoretic conditions on semi-Abelian varieties (Mumford-Tate, Kummer theory, Mordell-Weil, Schanuel Conjecture); we wanted an analysis covering more general algebraic varieties which would to lead to geometric conditions in place of those above.

3.9. **Thanks.** I thank my Mother and Father for support, patience and more. I also thank Artem Harmaty for attention to this work, and encouraging conversations, and Martin Bays for reading and discussing. Detailed thanks are in the report [Gavrilovich].

4. Appendix. Some examples.

Below we use cwf-notation and the language of model categories to give some examples in our model category $QtNaamen$ and the model category $Top$ of topological spaces. All claims we make are either standard or follow from the definitions and may be found in [Gavrilovich].

4.1. *Homotopy category. Cofibrant and fibrant objects.* Cofibrant objects, i.e. objects $X$ such that $\emptyset \xrightarrow{(c)} X$, are families of countable sets. A family $X$ is fibrant, i.e. $X \xrightarrow{(f)} \{x : x = x\}$, iff for every $x \in X$ and every $a$ finite, the union $\{x \cup a\}$ is also covered by a member of $X$, in notation $\{x \cup a\} \longrightarrow X$. We ignore non-existence of $\{x : x = x\}$ in ZFC; note that in ZFC fibrant objects are necessarily proper classes. The homotopy category $HoQtNaamen$ is, up to equivalence of categories: (i) $HoQtNaamen$ is the full subcategory of fibrant and cofibrant objects, (ii) $HoQtNaamen$ is the category of families of countable sets, with the arrows: $X \longrightarrow Y$ iff every $x \in X$ is almost covered by an element $y \in Y$, i.e. $X \dashrightarrow \cdot \xleftarrow{(wc)} Y$.

4.2. $StNaamen$ vs $QtNaamen$. Put a label $(q)$ on an arrow $A \longrightarrow B$ iff $A \longrightarrow B \curlywedge X \longrightarrow Y$ lifts wrt any arrow $X \longrightarrow Y$ between objects of $QtNaamen$. Then for any $A \in \mathcal{O}bStNaamen$ there exists an $\tilde{A} \in \mathcal{O}bQtNaamen$, unique up to isomorphism, such that $A \xrightarrow{(q)} \tilde{A}$. Diagram chasing using $(q)$-labels and M6 of StNaamen shows that the category $QtNaamen$ is closed under $M2$-decomposition, i.e. if $A, Y \in \mathcal{O}bQtNaamen$ and $A \xrightarrow{(wc)} B \xrightarrow{(f)} Y$ and $A \xrightarrow{(c)} X \xrightarrow{(wf)} Y$, then $B, X \in \mathcal{O}bQtNaamen$.

4.3. *Singletons.* For sets $A$ and $B$, $A \subseteq B$ iff there is a (necessarily unique) arrow/morphism $\{A\} \longrightarrow \{B\}$, and $\{A\} \xrightarrow{(wc)} \{B\}$ is an acyclic cofibration iff $B \setminus A$ is finite (and $A \subseteq B$). For sets $A$ and $B$ infinite $\operatorname{card} A = \operatorname{card} B$ iff $\{A\} \xrightarrow{(c)} \{B\}$, and $B$ is countable iff $\varnothing \xrightarrow{(c)} \{B\}$ is a cofibration.

4.4. *Ordinals.* For a ordinal it holds $\alpha \longrightarrow \{\alpha\}$ and $\cup \alpha \longrightarrow \alpha$.
$\{\alpha\} \longrightarrow \alpha + 1 \longrightarrow \{\alpha\}$, i.e. $\{\alpha\}$ and $\alpha + 1$ are isomorphic
$\alpha \xrightarrow{(f)} \alpha + 1$ iff $\alpha = \cup_{\beta < \alpha} \beta$ is limit.
$\alpha \xrightarrow{(wf)} \alpha + 1$ iff $\alpha$ is a regular cardinal, i.e. $\operatorname{cf}\alpha = \alpha$
$\alpha \xrightarrow{(c)} \alpha + 1$ iff $\alpha = \omega$ or $\alpha$ is not a cardinal
$\alpha \xrightarrow{(wc)} \alpha + 1$ iff $\alpha$ is not a limit ordinal, i.e. $\alpha \neq \cup_{\beta<\alpha}\beta$
$\alpha \xrightarrow{(c)} \beta$ iff $\alpha = \beta$ or $\alpha$ is not a cardinal and
  either $\operatorname{card}\beta \leqslant \operatorname{card}\alpha + \aleph_0$ or $\beta$ is a cardinal and $\operatorname{card}\beta \leqslant (\operatorname{card}\alpha + \aleph_0)^+$.
$\alpha \xrightarrow{(wc)} \beta$ iff $\beta < \alpha + \omega$ and $\alpha$ not a limit ordinal
$\alpha \in \mathcal{O}bQtNaamen$ iff $\operatorname{cf}\alpha = \omega$ or $\operatorname{cf}\alpha = \alpha$

4.5. *Fibrations. Increasing chains. Paths.* Take a set $M$ and represent $M$ as a union of a continuous increasing chain $M = \cup_{i<\lambda} M_i$; then $\{M_i\}_{i<\lambda} \xrightarrow{(f)} \{M\}$ is a fibration. Let $M^{<\lambda}$ be the set of all subsets of $M$ of cardinality strictly less than $\lambda$, then $M^{<\lambda} \xrightarrow{(f)} \{M\}$. If $\lambda > \omega$, then $M^{<\lambda} \xrightarrow{(wf)} \{M\}$ and $\{M_i\}_{i<\lambda} \xrightarrow{(wf)} \{M\}$ provided card $M_i <$ card $M$ and card $M = $ cf card $M$ is regular.

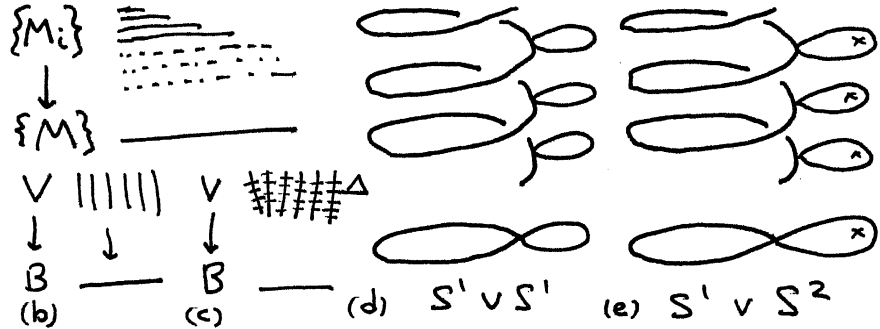

Fig. 1. Fibrations in $StNaamen$ and $Top$. (a) a union of increasing chain $\{M_i\}_{i<\lambda} \xrightarrow{(f)} \{M\}$ (b) a schematic picture of a fibration $V \xrightarrow{(f)} B$ of topological spaces; the vertical lines denote fibres $f^{-1}(b)$. (c) a schematic picture of a fibration $V \xrightarrow{(f)} B$ of topological spaces with a homotopy connexion $\Delta$, a rule for lifting uniquely paths $\gamma : [0,1]$ in $B$ to $\tilde\gamma : [0,1] \longrightarrow V$, $p(\tilde\gamma(t)) = \gamma(t)$ from an arbitrary point $y_0$, $\tilde\gamma(0) = y_0$. the horizontal lines represent paths so obtained. (d) a fibration whose base is a bouquet $S^1 \vee S^1$ of two circles (e) a fibration whose base is a bouquet $S^1 \vee S^2$ of a circle and a sphere.

4.6. *Cofibrations. Countable sets. Closed inclusions. Simplices.* We have $\emptyset \xrightarrow{(wc)} X$ is an acyclic cofibration iff $X$ is a family of finite sets, and $\emptyset \xrightarrow{(c)} X$ is a cofibration iff $X$ is a family of countable sets. An arrow $A \xrightarrow{(c)} B$ is a cofibration iff for every $b \in B$ there exists finitely many $b_0, ..., b_n = b \in B$, $n \in \omega$ and $a_0 \in A$ such that card $b \leqslant$ card $(a_0 \cap b_0)$ and card $b \leqslant$ card $(b_i \cap b_{i+1})$ for every $0 \leqslant i \leqslant n-1$. An arrow $A \xrightarrow{(wc)} B$ is an acyclic cofibration iff every $b \in B$ is almost a subset of an $a \in A$ (and of course $A \longrightarrow B$).

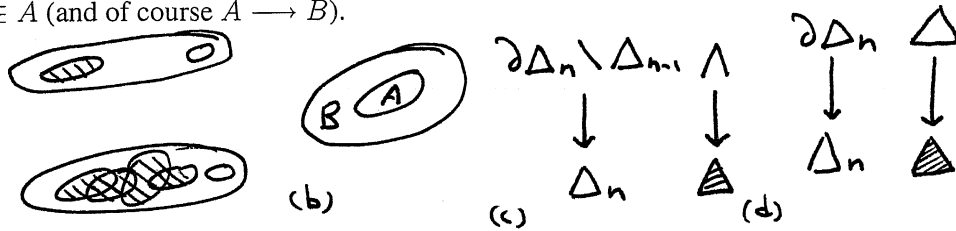

Fig. 2. Cofibrations in $StNaamen$ and $Top$. (a) a cofibration in $StNaamen$ (b) a picture of a cofibration $V \xrightarrow{(c)} B$ of topological spaces in Hurewicz model structure; this is just a closed inclusion. (c) a generating acyclic cofibration in Quillen model structure on Top: inclusion of the boundary without a face of a simplex into the simplex (d) a generating cofibration in Quillen model structure on Top: inclusion of the boundary of a simplex into the simplex

4.7. *Axiom M2. Path and cylinder spaces.* Let $A \longrightarrow Y$ be a morphism in $StNaamen$. The M2-decomposition can be explicitly given as follows:

$$A \xrightarrow{(wc)} \{ (a \cup y_{fini}) \cap y \ :\ a \in A,\ y \in Y,\ y_{fini} \text{ finite} \} \xrightarrow{(f)} Y$$

and

$$A \xrightarrow{(c)} \{ y \ :\ a_0 \in A,\ n \in \omega,\ y \subseteq y_n,\ y_0, ..., y_n \in Y,\ \operatorname{card} y \leqslant \operatorname{card}(a_0 \cap y_0) + \aleph_0$$

$$\text{and } \operatorname{card} y \leqslant \operatorname{card}(y_i \cap y_{i+1}) + \aleph_0 \text{ for every } 0 \leqslant i \leqslant n-1 \} \xrightarrow{(wf)} Y$$

In $Top$, for sufficiently nice topological spaces $A$ and $Y$, there are the following decompositions of a map $g: A \longrightarrow Y$:

$$A \xrightarrow{(wc)} \{ (a, \gamma) \ :\ a \in A,\ \gamma: [t_0, t_1] \longrightarrow Y,\ \gamma(t_0) = g(a) \} \xrightarrow{(f)} Y$$

$$A \xrightarrow{(c)} A \times [0,1] \cup Y/(a,0) \approx g(a) \xrightarrow{(wf)} Y$$

The maps involved are: in the wc-f-decomposition, a point $a \in A$ goes into the pair $(a, \gamma_{g(a)})$ where $\gamma$ is the constant path at point $g(a)$, $t_0 = t_1 = 0$; a pair $(a, \gamma)$ goes into $\gamma(t_2)$; and in the c-wf-decomposition, a point $a \in A$ goes into $(a, 1)$; and a point $(a, t) \in A \times [0, 1]$ goes into $g(a)$, and $y \in Y$ goes into itself.

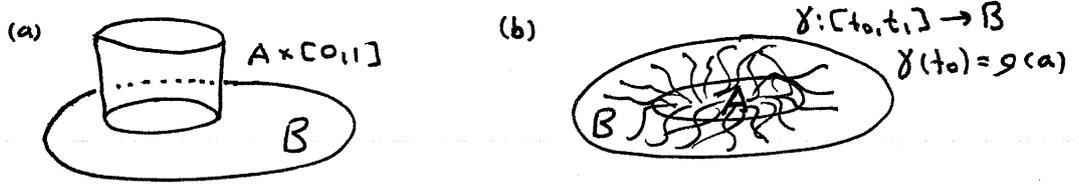

Fig. 3. M2-decompositions. Cones, cylinders and paths. (a) A c-wf-decomposition in Top using a cone object. (b) A wc-f-decomposition in Top using a cocone object of paths.

4.8. *Axiom M2. Downward Lowenheim-Skolem theorem as an instance of M2(c-wf).* Downward Lowenheim-Skolem theorem, e.g. for a first-order theory in a countable language, claims that every infinite subset of a model is contained in an elementary submodel of the same cardinality; equivalently, every subset $A \subseteq M$ is contained in an elementary submodel $A \subseteq M' \preceq M$ of $M$ of cardinality $\operatorname{card} M' = \operatorname{card} A + \aleph_0$. In our notation this is

$$\{A\} \xrightarrow{(c)} \{M' \ :\ M' \preceq M\} \xrightarrow{(wf)} \{M\}$$

In fact, it is easy to see that Downward Lowenheim-Skolem theorem is equivalent to the existence of the (c-wf)-decompositions in the category $QtNaamen(\mathfrak{R}_T)$ of §1.3.3 provided by Axiom M2 (where $\mathfrak{R}_T$ is the class of models of a countable language first-order theory $T$).

5. Appendix. Examples of lifting properties.

We give examples of some widely used notions that can be defined by a lifting property. Arguably, it is useful to think of these definitions as follows: we take a counterexample and "forbid" it by requiring the lifting property wrt to it. The following example may make this more clear. Assume we are interested in counting something, and we realise that to hope to preserve the count we need to avoid the two simplest possible(?) operations: adding a point $\{\} \longrightarrow \{\cdot\}$ to nothing or gluing two points into one $\{\cdot,\cdot\} \longrightarrow \{\cdot\}$. However, avoiding just these two is not enough: what we want is a class of operations(morphisms) which have nothing to do with these two bad ones. And we define such a class by requiring the left lifting property (Fig.4(b-c)). This gives us the class of bijections, i.e. exactly the operations that preserve the count.

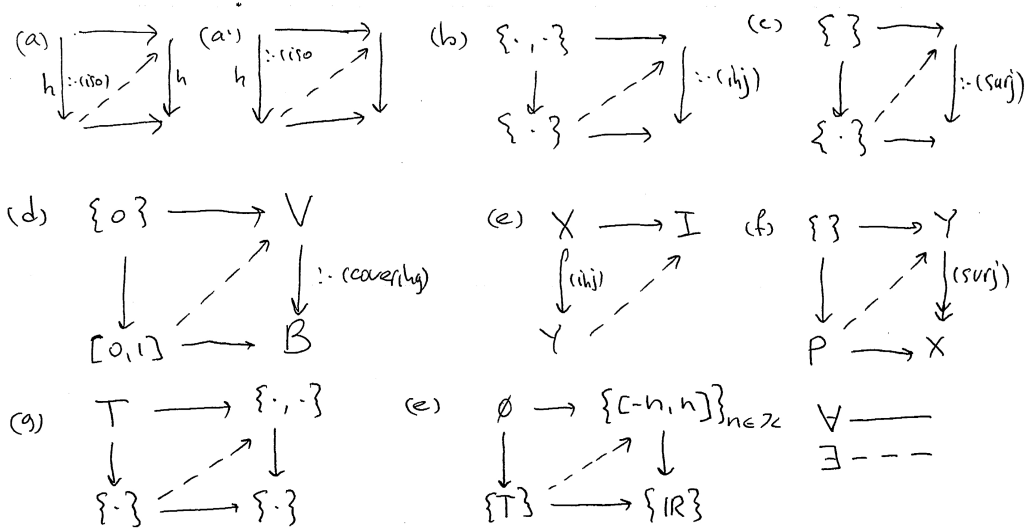

Fig.4. Read a diagramme as an $\forall\exists$-formula with parameters: for the arrows labelled ► or ◄, the following property holds: "for each commutative diagramme of solid arrows carrying labels as shown, there exists dashed arrows carrying labels as shown, making the diagramme of all the arrows commutative" (a) Isomorphism. In a category an arrow is an isomorphism iff it has (either left or right) lifting property wrt itself (and consequently ($a'$) any other arrow). (b) an arrow is injective iff it has the right lifting property wrt $\{\cdot,\cdot\}\longrightarrow\{\cdot\}$ whenever(=in most categories where) the latter notation/arrow makes sense. (c) an arrow is surjective iff it has left lifting property wrt $\{\}\longrightarrow\{\cdot\}$ whenever(=in most categories where) the latter notation/arrow makes sense. (d) Let $I = [0,1]$ be the unit interval of the real line, and let $0 \in [0,1]$ be its end point; the morphism $V\longrightarrow B$ is a *covering* of topological spaces iff there is always exists a unique lifting arrow $I\longrightarrow V$ making the diagramme commute. (e) an object $I$ is injective iff for each injective arrow $X\longrightarrow Y$ and any arrow $X\longrightarrow I$, there exists an arrow $Y\longrightarrow I$. (f) dually, an object $P$ is a projective object, e.g. a free module, iff for each surjective arrow $X\longleftarrow Y$ and an arrow $X\longleftarrow P$, there exists an arrow $Y\longleftarrow P$. (g) a topological space $T$ is connected iff $T\longrightarrow\{\cdot\}$ has the right lifting property wrt to $\{\cdot,\cdot\}\longrightarrow\{\cdot\}$ in the category of topological spaces) (e) a topological space $T$ is compact iff every continuous map $T\longrightarrow\mathbb{R}$ factors via an interval $[-n,n]$ for some $n \in \mathbb{Z}$.

# Axioms of a model category in labelled commutative diagrammes notation.

We state the axioms of Quillen of a model category in their original form. In particular, we follow the axiom numeration of Quillen(Homotopical Algebra).

**Notation** (Commutative diagrammes). *Commutative diagrammes will be used systematically throughout this note. Most importantly, diagrammes will be used to introduce new definitions. We introduce our notation for commutative diagrams. The properties defined are always properties of arrows. To distingish the arrows in the diagrammes which are the object of the definition we will denote them by ◄ or ►. We will mostly use commutative diagrammes to introduce ∀∃-definitions. In such cases solid arrows will be universally quantified and dashed arrows will be existentially quantified. Whenever definitions involving higher quantifier depth (such as in Figure ) a legend will be provided. As in Figure 1, we will use the notation $X \xrightarrow{\therefore (\cdot)} Y$ to mean "if the commutative diagram is true, then $X \longrightarrow Y$ is labeled $(\cdot)$". Notation $X \xrightarrow{!} Y$ indicates uniqueness. A legend on the right might be used to indicate the quantifiers and their order (from top to bottom). Unless stated otherwise, solid arrows are quantified universally, and dotted arrows are quantified existentially.*

**Axiom** ($M0$). *The category $\mathcal{C}$ is closed under finite projective and injective limits. It is known that it is enough to require existence of initial objects, terminal objects and pullbacks and pushouts.*

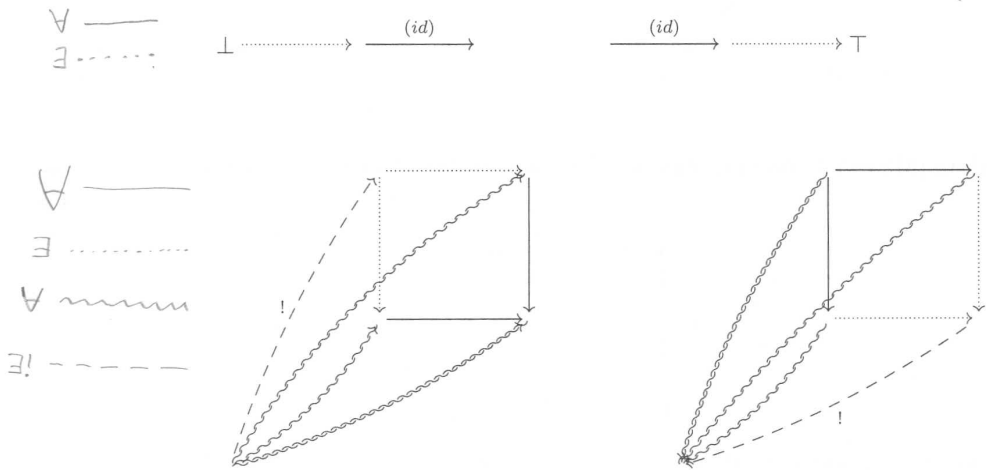

i

**Axiom** (M1). *The two following lifting properties for labeled arrows hold:*

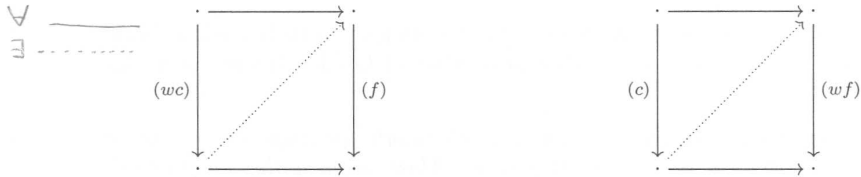

**Axiom** (M2). *The following two $\forall\exists$-diagrams hold:*

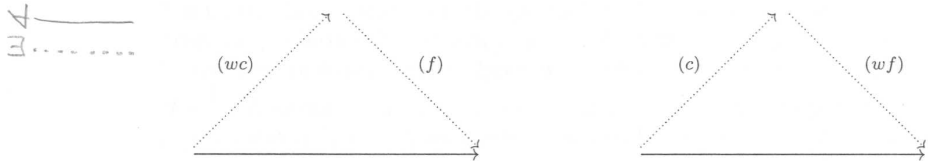

**Axiom** (M3(ccc,fff)). *Fibrations and cofibrations are stable under compositions. Namely, the following two $\forall\exists$-diagrams hold:*

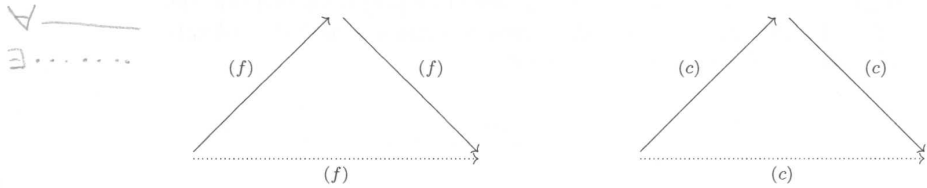

**Axiom** (M2(cwf)). *Isomorphisms are fibrations, co-fibrations and weak equivalences:*

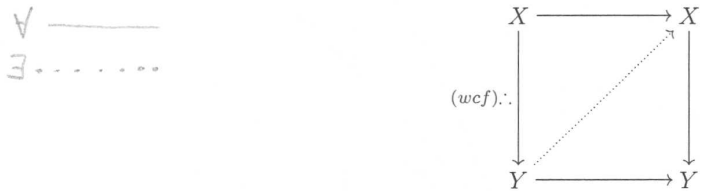

Figure 1: *The figure reads: if the commutative $\forall\exists$-diagramme is true then the left arrow is labeled (wcf).*

**Axiom** (M3($f \leftarrow f, c \rightarrow c$)). *Fibrations and cofibrations are stable under base change and co-base change respectively. I.e. the following diagrammes are true:*

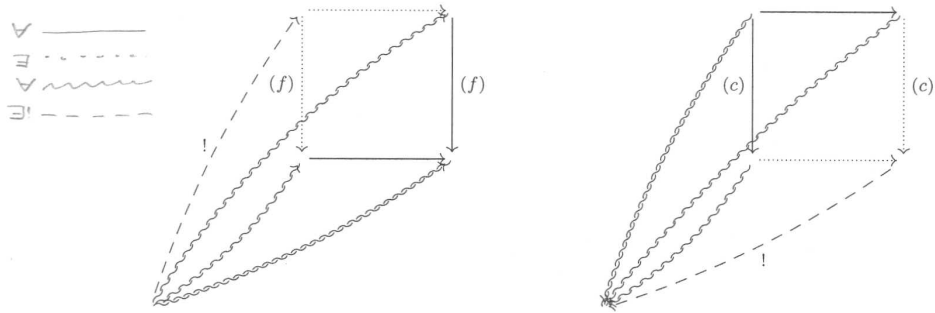

**Axiom** (M4($wf \leftarrow w, wc \rightarrow w$)). *The base extension of an arrow labeled (wc) and the co-base extension of an arrow labeled (wf) are both labeled (w):*

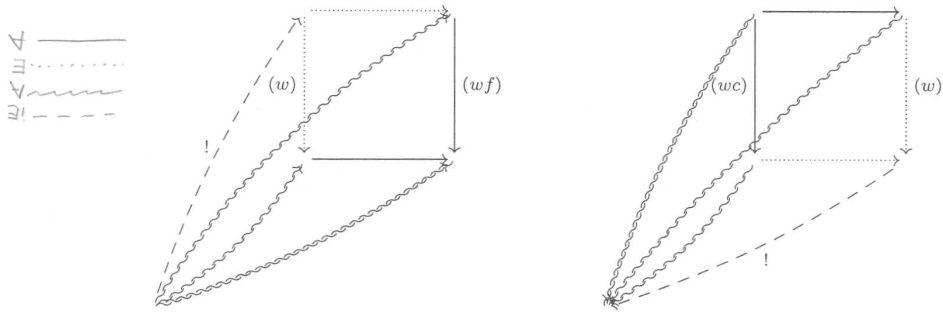

The last axiom assures that weak equivalence is close enough to being transitive:

**Axiom** (M5, Two out of three). *In a triangluar diagram, if any two of the arrows are labeled (w) so is the third*

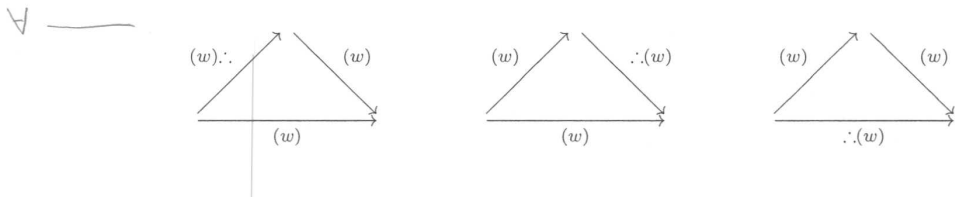